\documentclass[12pt]{article}

\usepackage{amsmath,amsthm}
\usepackage{amssymb}
\usepackage{bm}
\usepackage{graphicx}
\AtBeginDocument{
        \DeclareGraphicsExtensions{.pdf,.jpg,.png}}

\usepackage{enumitem}

\usepackage{hyperref}

\usepackage[T1]{fontenc}

\usepackage{listings}
\usepackage{color} 
\definecolor{mygreen}{RGB}{28,172,0} 
\definecolor{mylilas}{RGB}{170,55,241}

\newtheorem{theorem}{Theorem}[section]

\newtheorem{prop}[theorem]{Proposition}

\theoremstyle{definition}
\newtheorem{definition}[theorem]{Definition}
\newtheorem{remark}[theorem]{Remark}
\newtheorem{example}[theorem]{Numerical Example}
\newtheorem{conj}[theorem]{Conjecture}

\numberwithin{equation}{section}

\newcommand{\R}{\mathbb{R}}
\newcommand{\C}{\mathbb{C}}
\newcommand{\Z}{\mathbb{Z}}

\newcommand {\aplt} {\ {\raise-.5ex\hbox{$\buildrel<\over\sim$}}\ } 
\def\dddots{\mathinner{\mkern1mu\raise\p@  
    \hbox{.}\mkern2mu\raise4\p@\hbox{.}\mkern2mu
    \raise7\p@\vbox{\kern7\p@\hbox{.}}\mkern1mu}}%

\begin{document}


\baselineskip=17pt


\title{On Fekete Points for  a Real Simplex}

\author{Len Bos\\
Department of Computer Science\\ 
Univesity of Verona\\
Italy\\
\\
E-mail: leonardpeter.bos@univr.it}
\maketitle
\date{}

\renewcommand{\thefootnote}{}



\renewcommand{\thefootnote}{\arabic{footnote}}
\setcounter{footnote}{0}

\lstset{language=Matlab,%
    breaklines=true,%
    morekeywords={matlab2tikz},
    keywordstyle=\color{blue},%
    morekeywords=[2]{1}, keywordstyle=[2]{\color{black}},
    identifierstyle=\color{black},%
    stringstyle=\color{mylilas},
    commentstyle=\color{mygreen},%
    showstringspaces=false,
    numbers=left,%
    numberstyle={\tiny \color{black}},
    numbersep=9pt, 
    emph=[1]{for,end,break},emphstyle=[1]\color{red}, 
}

\begin{center}
Dedicated to Prof. Jaap Korevaar on the occasion of his 100th birthday!
\end{center}

\begin{abstract}
We survey what is known about Fekete points/optimal designs for a simplex in $\R^d.$ Several new results are included. The notion of Fej\'er exponenet for a set of interpolation points is introduced.
\end{abstract}

In one variable, as is well known, the Chebyshev points provide an excellent set of points for polynomial interpolation of functions defined on the interval $K=[-1,1].$ In several variables the problem of finding analogues of such near optimal interpolation points is much more difficult and each underlying set $K\subset \R^d$ must be analyzed indvidually.  One general approach, that turns out to be rather fruitful, is to consider the so-called Fekete points (see Definition \ref{DefFekete} below). These turn out to be strongly related to statistical Optimal Designs (cf. Definition \ref{DefOptDesign}) and to a property proved first by Fej\'er \cite{F32} (cf.  Defintion \ref{DefFejer})in the interval 
case. In this work we survey what is known about Fekete points for the case of $K\subset \R^d,$ a simplex.  Some new results are obtained and we introduce the notion of Fej\'er exponent for a set of interpolation points. We remark that also in the univariate complex case, $K\subset \C,$ Fekete points and their properties have been much studied, including by Prof.  J. Korevaar \cite{KM01}, (and the references therein),  to whom this paper is dedicated.

Consider then a set of $d+1$ points $X_d:=\{{\bf V}_1,\cdots,{\bf V}_{d+1}\}\subset \R^d$ in general position, and let $S_d:=\hbox{conv}(X_d)$ be the simplex generated from the vertices $X_d.$

We note that the dimension of the polynomials of degree at most $n$ in $d$  real variables is
\[{\rm dim}({\cal P}_n(\R^d))=N_n(=N):={n+d\choose d}.\]

For a basis $\{p_1,\cdots,p_N\}$ of  ${\cal P}_n(\R^d)$ and
$N$ points ${\bf x}_1,\cdots,{\bf x}_N$ in $S_d$ we may form the Vandermonde determinant
\[{\rm vdm}({\bf x}_1,\cdots,{\bf x}_N):={\rm det}([p_j({\bf x}_i)]_{1\le i,j\le N}).\]

In case the vandermonde determinant is non-zero,  the problem of interpolation at these points by polynomials of degree at most $n$ is regular,  and we may,  in particular,  construct the fundamental Lagrange polynomials $\ell_i({\bf x})$ of degree $n$ with the property that
\[\ell_i({\bf x}_j)=\delta_{ij}.\]

\begin{definition} \label{DefFekete} A set $F\subset S_d$ of $N$ distinct points is  said to a set of Fekete points of degree $n$ if they maximize 
$|{\rm vdm}({\bf x}_1,\cdots,{\bf x}_N)|$ over $S_d^N.$
\end{definition}

\begin{definition} \label{DefFejer} A set  $F\subset S_d$ of $N_d$ distinct points is said to be a Fej\`er set if
\[\max_{x\in S_d} \sum_{i=1}^{N_n}\ell_i^2(x)=1.\]
\end{definition}

More generally,  we may consider for a probability measure
\[\mu \in {\cal M}(S_d):=\{\mu\,:\, {\rm a\,\,probability\,\,measure\,\,supported\,\,on}\,\,S_d\},\]
the associated Gram matrix
\[ G_n(\mu):=[\int_{S_d}p_i(x)p_j(x)d\mu]\in \R^{N\times N}.\]

\begin{definition} \label{DefOptDesign} A measure $\mu\in {\cal M}(S_d)$ for which
${\rm det}(G_n(\mu)$ is a {\it maximum} is said to be D-optimal.
\end{definition}

By the Kiefer-Wolfowitz equivalence theorem \cite{KW} a measure is D-optimal if and only if it is G-optimal

\begin{definition}
A measure $\mu\in {\cal M}(S_d)$ for which the diagonal of the reproducing kernel
\[K_n({\bf x},{\bf x}):=\sum_{i=1}^N P_i^2({\bf x})\]
where $\{P_1,\cdots,P_N\}$ is any orthonormal basis for ${\cal P}_n(K)$ with respect to the inner product induced by $\mu,$ is  such that
\[ \max_{{\bf x}\in S_d}K_n({\bf x})\]
is a {\it minimum},  is said to be G-optimal.  In which case $\max_{{\bf x}\in S_d}K_n({\bf x})=N$ and this maximum is attained at all points in the support of $\mu.$
\end{definition}

For short,  we will refer to either a D-optimal or G-optimal measure as an {\it optimal} probability measure, or {\it optimal design},for  degree $n.$

The interested reader may find more about the theory of optimal designs in the monographs \cite{KS} and \cite{DS}. The asymptotics of such measures as $n\to\infty$ is discussed in \cite{BBLW}.

\begin{prop} Consider discrete probablity measures of the form
\[\mu =\sum_{{\bf a}\in F} w_{{\bf a}}\delta_{{\bf a}},  \quad w_{{\bf a}}>0\]
supported  on a set $F\subset S_d$ of cardinality $N.$ Then $\mu$ is an optimal measure for degree $n$ if and only if it is equallly weighted,  i.e.,  $w_{{\bf a}}=1/N,\,\,\forall \,{\bf a}\in F,$ and
$F$ is a Fej\`er set.
\end{prop}
\noindent {\bf Proof}.   For a set of cardinality $N$ it is easily seen that the associated Lagrange polynomials $\ell_{\bf a}({\bf x})/\sqrt{w_{{\bf a}}}$ are orthonormal with respect to the measure $\mu.$  Hence
\[K_n({\bf x})=\sum_{{\bf a}\in F}\frac{\ell_{\bf a}^2({\bf x})}{w_{\bf a}}.\]
Now suppose first that $\mu$ is optimal.  Then by G-optimality, for all ${\bf b}\in F,$
\begin{align*}
N&=K_n({\bf b})\cr
&=\sum_{{\bf a}\in F}\frac{\ell_{\bf a}^2({\bf b})}{w_{\bf a}}\cr
&=\frac{\ell_{\bf b}^2({\bf b})}{w_{\bf b}}\cr
&=\frac{1}{w_{\bf b}}
\end{align*}
and so the measure is equally weighted with $w_{\bf a}=1/N,$ $\forall {\bf a}\in F.$ Further,  by G-optimality,  we must also have $K_n({\bf x})\le N,$ ${\bf x}\in S_d.$ Hence
\begin{align*}
N&\ge K_n({\bf x})\cr
&=N \sum_{{\bf a}\in F}\ell_{\bf a}^2({\bf x}),
\end{align*}
i.e.,
\[ \sum_{{\bf a}\in F}\ell_{\bf a}^2({\bf x})\le 1,\quad {\bf x}\in S_d\]
and so $F$ is a Fej\`er set.

Conversely suppose that $F$ is a Fej\`er set and that  
\[ \mu =\frac{1}{N}\sum_{{\bf a}\in F}\delta_{{\bf a}}\] 
is equally weighted. Then 
\[K_n({\bf x})=N\sum_{{\bf a}\in F}\ell_{\bf a}^2({\bf x})\le N\times 1,\quad {\bf x}\in S_d\]
 and equal to $N$ at the points of $F.$ Hence the equally weighted measure $\mu$ is G-optimal.   $\square$

\begin{prop} A set of Fej\`er points $F$ is always also a set of Fekete points. 
\end{prop}
\noindent {\bf Proof}.   The equally weighted measure supported on $F,$
 \[ \mu =\frac{1}{N}\sum_{{\bf a}\in F}\delta_{{\bf a}}\]
 is by the previous proposition G-optimal,  and hence by the Kiefer-Wolfowitz equivalence theorem,  also D-optimal,  i.e.,  it maximizes the determinant of the Gram matrix for polynomials of degree at most $n$ among all probability measures.  But for an equally weighted discrete measure,  $\nu,$  supported on $N$ points ${\bf x}_1\cdots,{\bf x}_N,$  it is easy to confirm that
 \[ G_n(\nu)=\frac{1}{N}V_n^tV_n\]
 where $V_n:=[p_i(x_j)]_{1\le i\le n}$ is the Vandermonde matrix.
 Hence ${\rm det}(G_n(\nu))=(1/N^N)({\rm det}(V_n))^2$
 and so maximizing the Vandermonde determinant is equivalent to maximizing the Gram matrix over all {\it equally weighted} measures supported on $N$ points in $S_d.$ Since the measure supported on a set of Fej\`er points is optimal among {\it all} probability measures,  it is also best among equally weighted measures,  and hence $F$ is also a set of Fekete points.  $\square$
 
 \medskip

\begin{remark} It is not in general true that Fekete points by themselves are always Fej\`er points; see \cite{Bo} for a discussion of this problem.  $\square$ \end{remark}

\medskip\begin{remark}  Since a Lagrange polynomial may be written as a ratio of Vandermonde determinants
\[\ell_i({\bf x})=\frac{{\rm vdm}({\bf x}_1,\cdots,{\bf x}_{i-1},{\bf x},
{\bf x}_{i+1},\cdots,{\bf x}_N)}{{\rm vdm}({\bf x}_1,{\bf x}_{2}\cdots,{\bf x}_N)}\]
it is necessarily the case that,  for Fekete points,  the Lagrange polynomials are bounded by 1 on $S_d$ or, in other words,
\[\|{\bm \ell}({\bf x})\|_\infty\le 1,\quad {\bf x}\in S_d\]
where 
\[{\bm \ell}({\bf x})=[\ell_1({\bf x}),\cdots,\ell_N({\bf x})]\in\R^N\]
denotes the vector of the Lagrange polynomials.  Note however that this boundedness condition is {\it not} sufficient to be a Fekete set.  Indeed,  consider for degree 1 the simplex (triangle) $S_2$  in $\R^2,$ degree 1 and the three edge midpoints,  $(V_1+V_2)/2,$ $(V_1+V_3)/2$ and $(V_2+V_3)/2.$ It is easy to confirm that the associated Lagrange polynomials are $1-2\lambda_3,$ $1-2\lambda_2$ and $1-2\lambda_1,$ respectively.  As $\lambda_j\in[0,1],$ each of these Lagrange polynomials takes values in $[-1,1],$ i.e.,  is bounded by 1 in absolute value.  However,  the three edge midpoints do {\it not} form a Fekete set,  as their Vandermonde determinant is $1/2\times {\rm area}(S_2)$ whereas the Vandermonde determinant for the three vertices is $2\times {\rm area}(S_2).$ $\square$ \end{remark}

\begin{center}\begin{figure}[t!]\centering
 \includegraphics[width=.6\textwidth,height=.3\textheight]{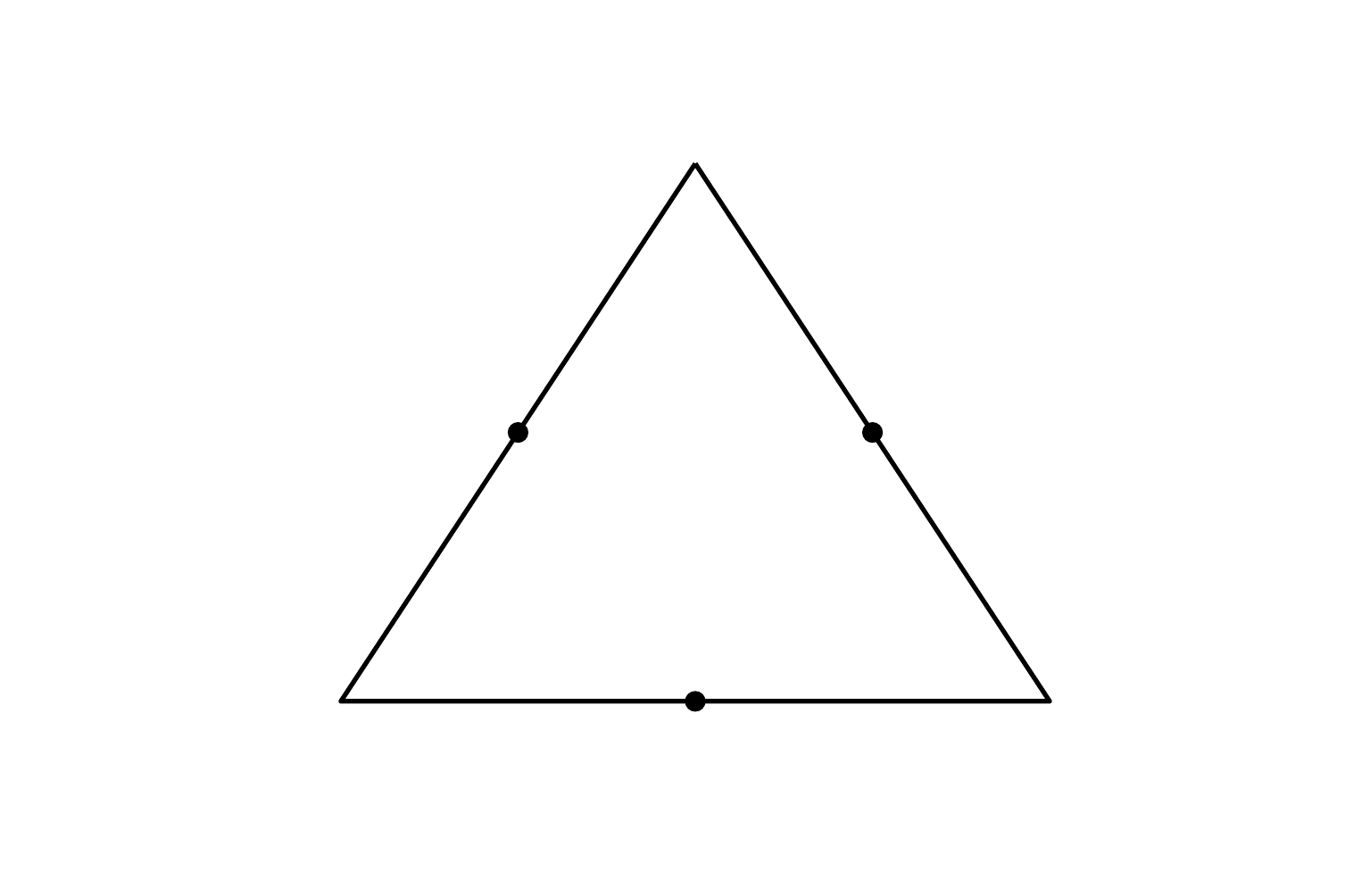}\hfill
 \caption{Three Edge Midpoints}
\end{figure}
\end{center}

\medskip We will now proceed to discuss the degrees 1 through 5 situations,  based on sets of points introduced in \cite{B83}.  The idea is to place a maximal number of points on the boundary of the simplex in $\R^d$,  in such a way that restricted to a lower  dimensional face of dimension $0\le d'<d$ we obtain the prescribed set of points for the simplex in $\R^{d'}.$ On any  {\it edge},  i.e. , a one dimensional face,  we place the univariate Fekete points,  which are known to be the zeros of $(x^2-1)P_n'(x),$ where $P_n(x)$ is the $n$th classical Legendre polynomial. 

\medskip

\section{Degree One}
For degree one the dimension of the polynomials is $d+1.$
\begin{prop} (cf. \cite{B83})
The set $F_1=X_d\subset S_d$ is a set of Fej\`er points of degree one for all dimensions $d.$ 
\end{prop}
\noindent {\bf Proof.} The Lagrange polynomial for the $i$th vertex $V_i$ is just
$\ell_i=\lambda_i,$ in barycentric coordinates. 
Now,  for $x\in S_d,$ $\lambda_i\ge0,$ $1\le i\le d+1$ and hence 
\begin{align*}
\sum_{i=1}^{d+1}\ell_i^2&=\sum_{i=1}^{d+1}\lambda_i^2\\
&\le \left(\sum_{i=1}^{d+1}\lambda_i \right)^2\\
&=1.
\end{align*}

$\square$

\section{Degree Two}
Let $F_2\subset S_d$ be the set consisting of the $d+1$ vertices together with the 
${d+1 \choose 2}$ edge midpoints. Note that 
\[(d+1)+ {d+1 \choose 2} = {d+2 \choose 2}=N_2.\]
\begin{prop} (cf. \cite{B83})
The set $F_2$ is a Fej\`er set for all dimensions $d.$
\end{prop}
\noindent {\bf Proof.} It is easy to confirm that the Lagrange polynomial for vertex $V_i$ is $\lambda_i(2\lambda_i-1)$ while for the midpoint $(V_i+V_j)/2$ it is 
$4\lambda_i\lambda_j.$ 

We make use of the symmetric power functions
\[p_r({\bm \lambda}):=\sum_{i=1}^{d+1}\lambda_i^r.\]
Note that for barycentric coordinates $p_1({\bm \lambda})=1,$ and hence
\begin{align*}
1-K_2({\bm \lambda})&=p_1^4-\sum_{i=1}^{d+1} \lambda_i^2(2\lambda_i-p_1)^2-16\sum_{1\le i<j\le d+1} \lambda_i^2\lambda_j^2\cr
&=p_1^4-\sum_{i=1}^{d+1}\bigl\{4\lambda_i^4-4\lambda_i^3 p_1 +\lambda_i^2p_1^2\bigr\} -8(p_2^2-p_4)\cr
&=p_1^4-\bigl\{4p_4-4p_1p_3+p_1^2p_2\bigr\}-8(p_2^2-p_4)\cr
&=p_1^4+4p_4+4p_1p_3-p_1^2p_2-8p_2^2.
\end{align*}
It can be shown,  by brute force calculus,  that this is always positive on the simplex $S_d.$ 

However,  a perhaps more elegant solution is to notice,  after some computation,   that
\begin{align*}1-K_2({\bm \lambda})&=6\sum_{j<k}\lambda_j\lambda_k
(\lambda_j-\lambda_k)^2 +10\sum_{i<j<k}\lambda_i\lambda_j\lambda_k(\lambda_i+\lambda_j+\lambda_k) \\
&+24\sum_{i<j<k<m}\lambda_i\lambda_j\lambda_k\lambda_m\\
&\ge0,\quad {\bf x}\in S_d.
\end{align*}
$\square$

\begin{remark} It is worthwhile noting that $1-K_2({\bf x})$ can also be analyzed using Schur functions.  Indeed,  
for a partition ${\bm \mu}\in \Z_+^m$ with
\[\mu_1\ge\mu_2\ge\cdots\ge \mu_m\ge0\]
one defines (cf.  the monograph by Macdonald \cite{M95}),  for ${\bf x}\in \R^m,$
\begin{equation}\label{slambda}
s_{\bm \mu}({\bf x})={\rm det}([x_i^{\mu_j+m-j}]_{1\le i,j\le m})/\prod_{1\le i<j\le m}(x_i-x_j).
\end{equation}
The Schur functions are symmetric polynomials and indeed the Schur polynomials $s_{\bm \mu}(x_1,\cdots,x_m)$ with $|{\bm \mu}|=n,,$ form a $\Z$-basis for the homogeneous symmetric polynomials of degree $n$ with integer coefficients in $m$ variables.

In particular we may express $1-K_2({\bm \lambda})$ in terms of
Schur polynomials,  and indeed one may verify that in dimension $d$ with ${\bm \lambda} \in \R^{d+1},$
\begin{equation}\label{K2Schur}
1-K_2=6s_{[3,1,0^{d-1}]}-18 s_{[2,2,0^{d-1}]}+16 s_{[2,1,1,0^{d-2}]}-6s_{[1,1,1,1,0^{d-3}]}.
\end{equation}

Now it can also be verified that
\[\sum_{i<j} \lambda_i\lambda_j(\lambda_i-\lambda_j)^2 = s_{[3,1,0^{d-1}]}-3s_{[2,2,0^{d-1}]}+s_{[2,1,1,0^{d-2}]}\]
and so,  in particular,
\[ s_{[3,1,0^{d-1}]}-3s_{[2,2,0^{d-1}]}+s_{[2,1,1,0^{d-2}]}\ge0.\]
Hence by (\ref{K2Schur}),
\begin{align*}
1-K_2&=6( s_{[3,1,0^{d-1}]}-3s_{[2,2,0^{d-1}]}+s_{[2,1,1,0^{d-2}]})+10 s_{[2,1,1,0^{d-2}]}-6s_{[1,1,1,1,0^{d-3}]}\cr
&\ge 10 s_{[2,1,1,0^{d-2}]}-6s_{[1,1,1,1,0^{d-3}]}.
\end{align*}

We claim that the last expression is non-negative.  To see this we use the inequality on normalized Schur functions given by Sra \cite{S16}.  Specifically,  write 
\[{\bm \mu}_1\succeq {\bm \mu}_2 \,\,\iff\,\, \sum_{j=1}^k({\bm \mu_1})_j\ge \sum_{j=1}^k ({\bm \mu}_2)_k,\,\,k=1,2\cdots.\]
Then, in case ${\bm \mu}_1\succeq {\bm \mu}_2,$ for ${\bm x}\in \R_+^m,$
\begin{equation}\label{sra}
\frac{s_{{\bm \mu}_1}({\bf x})}{s_{{\bm \mu}_1}(1^m)} \ge
\frac{s_{{\bm \mu}_2}({\bf x})}{s_{{\bm \mu}_2}(1^m)}.
\end{equation}
We may compute
\[ s_{[2,1,1,0^{d-2}]}(1^{d+1})=3 {d+2 \choose 4},\,\,\, s_{[1,1,1,1,0^{d-3}]}(1^{d+1})={d+1 \choose 4}.\]

Consequently,  
\[s_{[2,1,1,0^{d-2}]}\ge 3\frac{d+2}{d-2} \times s_{[1,1,1,1,0^{d-3}]}\]
from which our claim follows easily.  
$\square$
\end{remark}

\section{Degree Three}
Let $F_3\subset S_d$ be the set consisting of the $d+1$ vertices together with the
$2{d+1 \choose 2}$  edge points  $tV_i+(1-t)V_j,$ $i\neq j,$ $t=(1+1/\sqrt{5})/2,$ together  with the ${d+1 \choose 3}$ barycentres of each two-dimensional face.

Note that
\[(d+1)+2{d+1 \choose 2}+{d+1 \choose 3}={d+3 \choose 3}=N_3.\]

We conjecture that $F_3$ is a  Fej\`er set for all dimensions, but can only prove this for dimensions up to 28.

\begin{prop} (cf. \cite{B83}) \label{deg3case} $F_3$ is a Fej\`er set for dimensions $1\le d\le 28.$
\end{prop}
\noindent {\bf Proof.} The $d=1$ case is a special case of Fej\`er's original theorem \cite{F32}.  Otherwise,  we note (as given in \cite{B83}) that for the vertex $V_i,$
\[\ell_i=(\lambda_i/2)(12\lambda_i^2-12\lambda_i+3-S_2)\]
where $S_k:=\sum_{i=1}^{d+1}\lambda_i^k;$ while for the face midpoint
$(V_i+V_j+V_k)/3,$ $i<j<k,$
\[\ell_{ijk}=27\lambda_i\lambda_j\lambda_k,\]
and for the edge points $tV_i+(1-t)V_j,$ $i\neq j,$
\[\ell_{ij}=5\lambda_i\lambda_j((1+t\sqrt{5})\lambda_i+(2-t\sqrt{5})\lambda_j-1).\]
The  $d=2$ case is special,  as there is a single point,  in the interior of $S_2$ (at the barycentre).  We let, for simplicity's sake,  $x,y,z$ denote  the three barycentric coordinates and compute
\begin{align*}
K(x,y,z)&:=\sum_{i=1}^{10} \ell_i^2(x,y,z)\\
&={z}^{6}-6\,y\,{z}^{5}-6\,x\,{z}^{5}+87\,{y}^{2}\,{z}^{4}+24\,x\,y\,{z}^{4}+87\,{x}^{2}\,{z}^{4}-112\,{y}^{3}\,{z}^{3}\\
&-68\,x\,{y}^{2}\,{z}^{3}-68\,{x}^{2}\,y\,{z}^{3}-112\,{x}^{3}\,{z}^{3}+87\,{y}^{4}\,{z}^{2}-68\,x\,{y}^{3}\,{z}^{2}\\
&+912\,{x}^{2}\,{y}^{2}\,{z}^{2}-68\,{x}^{3}\,y\,{z}^{2}+87\,{x}^{4}\,{z}^{2}-6\,{y}^{5}\,z+24\,x\,{y}^{4}\,z-68\,{x}^{2}\,{y}^{3}\,z\\
&-68\,{x}^{3}\,{y}^{2}\,z+24\,{x}^{4}\,y\,z-6\,{x}^{5}\,z+{y}^{6}-6\,x\,{y}^{5}+87\,{x}^{2}\,{y}^{4}-112\,{x}^{3}\,{y}^{3}\\
&+87\,{x}^{4}\,{y}^{2}-6\,{x}^{5}\,y+{x}^{6}
\end{align*}
so that 
\begin{align*}
H(x,y,z)&:=1-K(x,y,z)\\
&=(x+y+z)^6-K(x,y,z)\\
&=12\,y\,{z}^{5}+12\,x\,{z}^{5}-72\,{y}^{2}\,{z}^{4}+6\,x\,y\,{z}^{4}-72\,{x}^{2}\,{z}^{4}+132\,{y}^{3}\,{z}^{3}\\
&+128\,x\,{y}^{2}\,{z}^{3}+128\,{x}^{2}\,y\,{z}^{3}+132\,{x}^{3}\,{z}^{3}-72\,{y}^{4}\,{z}^{2}+128\,x\,{y}^{3}\,{z}^{2}\\
&-822\,{x}^{2}\,{y}^{2}\,{z}^{2}+128\,{x}^{3}\,y\,{z}^{2}-72\,{x}^{4}\,{z}^{2}+12\,{y}^{5}\,z+6\,x\,{y}^{4}\,z\\
&+128\,{x}^{2}\,{y}^{3}\,z+128\,{x}^{3}\,{y}^{2}\,z+6\,{x}^{4}\,y\,z+12\,{x}^{5}\,z+12\,x\,{y}^{5}\\
&-72\,{x}^{2}\,{y}^{4}+132\,{x}^{3}\,{y}^{3}-72\,{x}^{4}\,{y}^{2}+12\,{x}^{5}\,y.
\end{align*}
We are claiming that $H(x,y,z)\ge 0$ for $x,y,z\ge0,$ $x+y+z=1.$
Certainly this would be the case if all the coefficients were non-negative.  However,  $H(x,y,z)=0$ at the barycentre and hence it is not possible to have all non-negative coefficients.  To overcome this problem we re-express $H$ in terms of the barycentric coordinates for the subtriangle with vertices $V_1,$ $V_2$ and the barycentre $(V_1+V_2+V_3)/3.$ If $H\ge0$ on this subtriangle then, by symmetry,  it will be non-negative on the whole simplex.  

Let therefore $(u,v,w)$ be the barycentric coordinates with repsect to this subtriangle. Then
\[x=u+w/3,\,\,y=v+w/3,\,\,z=w/3\]
and $H$ becomes
\begin{align*}
H(u,v,w)&=\frac{2}{81}\{ 131\,{v}^{2}\,{w}^{4}-131\,u\,v\,{w}^{4}+131\,{u}^{2}\,{w}^{4}+312\,{v}^{3}\,{w}^{3}+318\,u\,{v}^{2}\,{w}^{3}\\
&+318\,{u}^{2}\,v\,{w}^{3}+312\,{u}^{3}\,{w}^{3}-81\,{v}^{4}\,{w}^{2}+1566\,u\,{v}^{3}\,{w}^{2}+1215\,{u}^{2}\,{v}^{2}\,{w}^{2}\\
&+1566\,{u}^{3}\,v\,{w}^{2}-81\,{u}^{4}\,{w}^{2}+324\,{v}^{5}\,w-1053\,u\,{v}^{4}\,w+3186\,{u}^{2}\,{v}^{3}\,w\\
&+3186\,{u}^{3}\,{v}^{2}\,w-1053\,{u}^{4}\,v\,w+324\,{u}^{5}\,w+486\,u\,{v}^{5}-2916\,{u}^{2}\,{v}^{4}\\
&+5346\,{u}^{3}\,{v}^{3}-2916\,{u}^{4}\,{v}^{2}+486\,{u}^{5}\,v\}.
\end{align*}
Then 
\begin{align*}
H=&\frac{262\,\left( {v}^{2}-u\,v+{u}^{2}\right) \,{w}^{4}}{81}+12\,u\,v\,\left( {v}^{2}-3\,u\,v+{u}^{2}\right)^{2}+\frac {2w}{27}
\{ 104\,{v}^{3}\,{w}^{2}\\
&+106\,u\,{v}^{2}\,{w}^{2}+106\,{u}^{2}\,v\,{w}^{2}+104\,{u}^{3}\,{w}^{2}-27\,{v}^{4}\,w+522\,u\,{v}^{3}\,w\\
&+405\,{u}^{2}\,{v}^{2}\,w
+522\,{u}^{3}\,v\,w-27\,{u}^{4}\,w+108\,{v}^{5}\\
&-351\,u\,{v}^{4}+1062\,{u}^{2}\,{v}^{3}+1062\,{u}^{3}\,{v}^{2}-351\,{u}^{4}\,v+108\,{u}^{5}\}.
\end{align*}
To show that $H\ge0$ for $u,v,w\ge0,$ $u+v+w=1,$ we note that as the first two terms above are both positive, it suffices to show that the expression in the brace brackets is also positive. Let $E$ denote this expression.

However,
\begin{align*}
(u+v+w)^4E&=104\,{v}^{3}\,{w}^{6}+106\,u\,{v}^{2}\,{w}^{6}+106\,{u}^{2}\,v\,{w}^{6}+104\,{u}^{3}\,{w}^{6}\\
&+389\,{v}^{4}\,{w}^{5}
+1362\,u\,{v}^{3}\,{w}^{5}+1253\,{u}^{2}\,{v}^{2}\,{w}^{5}+1362\,{u}^{3}\,v\,{w}^{5}\\
&+389\,{u}^{4}\,{w}^{5}+624\,{v}^{5}\,{w}^{4}+3513\,u\,{v}^{4}\,{w}^{4}\\
&+7302\,{u}^{2}\,{v}^{3}\,{w}^{4}+7302\,{u}^{3}\,{v}^{2}\,{w}^{4}+3513\,{u}^{4}\,v\,{w}^{4}+624\,{u}^{5}\,{w}^{4}\\
&+686\,{v}^{6}\,{w}^{3}+3508\,u\,{v}^{5}\,{w}^{3}+14320\,{u}^{2}\,{v}^{4}\,{w}^{3}+22996\,{u}^{3}\,{v}^{3}\,{w}^{3}\\
&+14320\,{u}^{4}\,{v}^{2}\,{w}^{3}+3508\,{u}^{5}\,v\,{w}^{3}+686\,{u}^{6}\,{w}^{3}+644\,{v}^{7}\,{w}^{2}\\
&+1476\,u\,{v}^{6}\,{w}^{2}+11522\,{u}^{2}\,{v}^{5}\,{w}^{2}
+31694\,{u}^{3}\,{v}^{4}\,{w}^{2}\\
&+31694\,{u}^{4}\,{v}^{3}\,{w}^{2}+11522\,{u}^{5}\,{v}^{2}\,{w}^{2}+1476\,{u}^{6}\,v\,{w}^{2}+644\,{u}^{7}\,{w}^{2}\\
&+405\,{v}^{8}\,w+306\,u\,{v}^{7}\,w+3663\,{u}^{2}\,{v}^{6}\,w+18378\,{u}^{3}\,{v}^{5}\,w\\
&+29232\,{u}^{4}\,{v}^{4}\,w+18378\,{u}^{5}\,{v}^{3}\,w+3663\,{u}^{6}\,{v}^{2}\,w+306\,{u}^{7}\,v\,w\\
&+405\,{u}^{8}\,w+108\,{v}^{9}+81\,u\,{v}^{8}+306\,{u}^{2}\,{v}^{7}+3636\,{u}^{3}\,{v}^{6}+8973\,{u}^{4}\,{v}^{5}\\
&+8973\,{u}^{5}\,{v}^{4}+3636\,{u}^{6}\,{v}^{3}+306\,{u}^{7}\,{v}^{2}+81\,{u}^{8}\,v+108\,{u}^{9}
\end{align*}
which is positive as each term has a positive coefficient. 

For $d\ge3$ there are no longer interior points and there is a simple algebraic procedure to verify that $1-K_3\ge 0$ on the simplex.   We proceed by induction.  Assuming that $F_3$ is a  Fej\`er set for dimension $d-1$ one calculates $1-K_3$ in homogenized form 
\[\Bigl(\sum_{j=1}^{d+1}\lambda_j\Bigr)^6-K_3({\bm \lambda})\]
and then considers this on the sub-simplex with vertices 
\[V_1,\cdots,V_d,  (V_1+\cdots+V_{d+1})/(d+1).\]
If we let $\mu_1,\cdots,\mu_{d+1}$ be the barycentric coordinates for the sub-simplex,  then it is easy to see that
\begin{equation}\label{mus}
\lambda_j=\mu_j+\mu_{d+1}/(d+1), \,\,1\le j\le d, \quad \lambda_{d+1}=\mu_{d+1}/(d+1).
\end{equation}
Let $H({\bm \mu})$ be $1-K_3$ restricted to this sub-simplex,  i.e.,
\[H({\bm \mu})=\Bigl(\sum_{j=1}^{d+1}\lambda_j\Bigr)^6-K_3({\bm \lambda})\]
with ${\bm \lambda}$ given by (\ref{mus}).  It is sufficient to prove that $H({\bm \mu})\ge 0$ for $\mu_j\ge0,$ $1\le j\le (d+1)$ for which it is, in turn, sufficient to show that
\[H({\bm \mu})-H(\widehat{{\bm \mu}},0)\ge 0,\quad \widehat{{\bm \mu}}:=(\mu_1,\mu_2,\cdots,\mu_d)\]
as the term $H(\widehat{{\bm \mu}},0)$ is the restriction of $H$ to the face $\mu_{d+1}=\lambda_{d+1}=0$ and hence is positive by the induction hypothesis. 

For positivity it is sufficient that, for some integer $r\ge0,$ all the coefficients of
\[ \Bigl(\sum_{j=1}^{d+1}\mu_j\Bigr)^r(H({\bm \mu})-H(\widehat{{\bm \mu}},0))\]
are non-negative.  This has been verified by means of computer algebra for degrees up to 28 with $r=8$ for dimension $d=3,$ $r=5$ for $d=4,$ $r=3$ for $d=5,6$ and $r=2$ otherwise.  Appendix A gives Matlab code using its Symbolic Toolbox for this purpose.  $\square$

\begin{remark} Multiplying a polynomial in barycentric coordinates by one in the form of $\sum_{j=1}^{d+1}\mu_j$ is known as {\it degree elevation}.  The minimal $r$ necessary to have all non-negative coefficients is connected to the so-called Bernstein degree of a polynomial which is positive on a simplex.  $\square$
\end{remark}


\section{Degree Four}

Let $F_4\subset S_d$ be the set of points introduced in \cite{B83} consisting of
\begin{itemize}
\item the ${d+1 \choose 4}$ centres of each  three dimesnional face: \[\frac{1}{4}(V_i+V_j+V_k+V_\ell),\,\, i<j<k<\ell\]
\item the $3{d+1 \choose  3}$ vertices of a triangle in the interior of each two dimensional face:
\[\frac{4+\sqrt{5}}{11}V_i+\frac{7-\sqrt{5}}{22}(V_j+V_k),\,\, i\neq j,k,\,\,j<k\]
\item the $3{d+1 \choose 2}$ points with 3 points on each edge (1 dimensional face):
\begin{align*}
&\frac{1-\sqrt{3/7}}{2}V_i+\frac{1+\sqrt{3/7}}{2}V_j,\,\, i\neq j,\cr
&\frac{1}{2}(V_i+V_j),\,\,i<j
\end{align*}
\item the $d+1$ vertices $V_i.$
\end{itemize}
We note that the cardinality of $F_4$ is
\[ {d+1 \choose 4}+3{d+1 \choose 3}+3{d+1\choose 2}+(d+1)={d+4 \choose 4}=N_4.\]
The associated Lagrange polynomials are also given in \cite{B83}, which,  for the sake of completeness,  we reproduce here:
\begin{itemize}
\item $\ell_{ijk\ell}=256\lambda_i\lambda_j\lambda_k\lambda_\ell$
\item $\displaystyle{\ell_{ijk}=\frac{73+\sqrt{5}}{2}\lambda_i\lambda_j\lambda_k \Bigl\{(1+\sqrt{5})\lambda_i+\frac{3-\sqrt{5}}{2}(\lambda_j+\lambda_k)-1\Bigr\}}$
\item \begin{align*}
\ell_{ij}&= \frac{49}{6}\lambda_i\lambda_j \Bigl\{ \bigl(
\frac{61+7\sqrt{5}}{22}-\sqrt{3/7}\frac{7+\sqrt{5}}{2}\bigr)\lambda_i^2 \cr
&+\frac{3+9\sqrt{5}}{11}\lambda_i\lambda_j + \bigl(\frac{61+7\sqrt{5}}{22}+\sqrt{3/7}\frac{7+\sqrt{5}}{2}\bigr)\lambda_j^2\cr
&+\bigl(\sqrt{3/7}\frac{5+\sqrt{5}}{2}-\frac{41+13\sqrt{5}}{22}\bigr)\lambda_i\cr
&-\bigl(\sqrt{3/7}\frac{5+\sqrt{5}}{2}+\frac{41+13\sqrt{5}}{22}\bigr)\lambda_i\cr
&+\frac{8+2\sqrt{5}}{11}-\frac{7-\sqrt{5}}{11}S_2\Bigl\}
\end{align*}
\item \begin{align*}
\ell_{ij}&= \frac{4}{33}\lambda_i\lambda_j \Bigl\{ (100-8\sqrt{5})(\lambda_i^2+\lambda_j^2)+(672-96\sqrt{5})\lambda_i\lambda_j\cr
&+(62\sqrt{5}-302)(\lambda_i+\lambda_j)+(82-40\sqrt{5})S_2\Bigr\}
\end{align*}
\item \begin{align*}
\ell_i&=2\lambda_i\Bigl\{\frac{101+17\sqrt{5}}{11}\lambda_i^3
-(13+2\sqrt{5})\lambda_i^2+\frac{130+27\sqrt{5}}{22}\lambda_i\cr
&-\frac{259+51\sqrt{5}}{264} +S_2\bigl( \frac{9+\sqrt{5}}{8}-\frac{24+17\sqrt{5}}{22}\lambda_i\bigr)\cr
&+S_3\bigl(\frac{45+3\sqrt{5}}{44}-\frac{5}{3}\bigr)\Bigr\}.
\end{align*}
\end{itemize}

\begin{center}\begin{figure}[t!]\centering
 \includegraphics[width=.6\textwidth,height=.3\textheight]{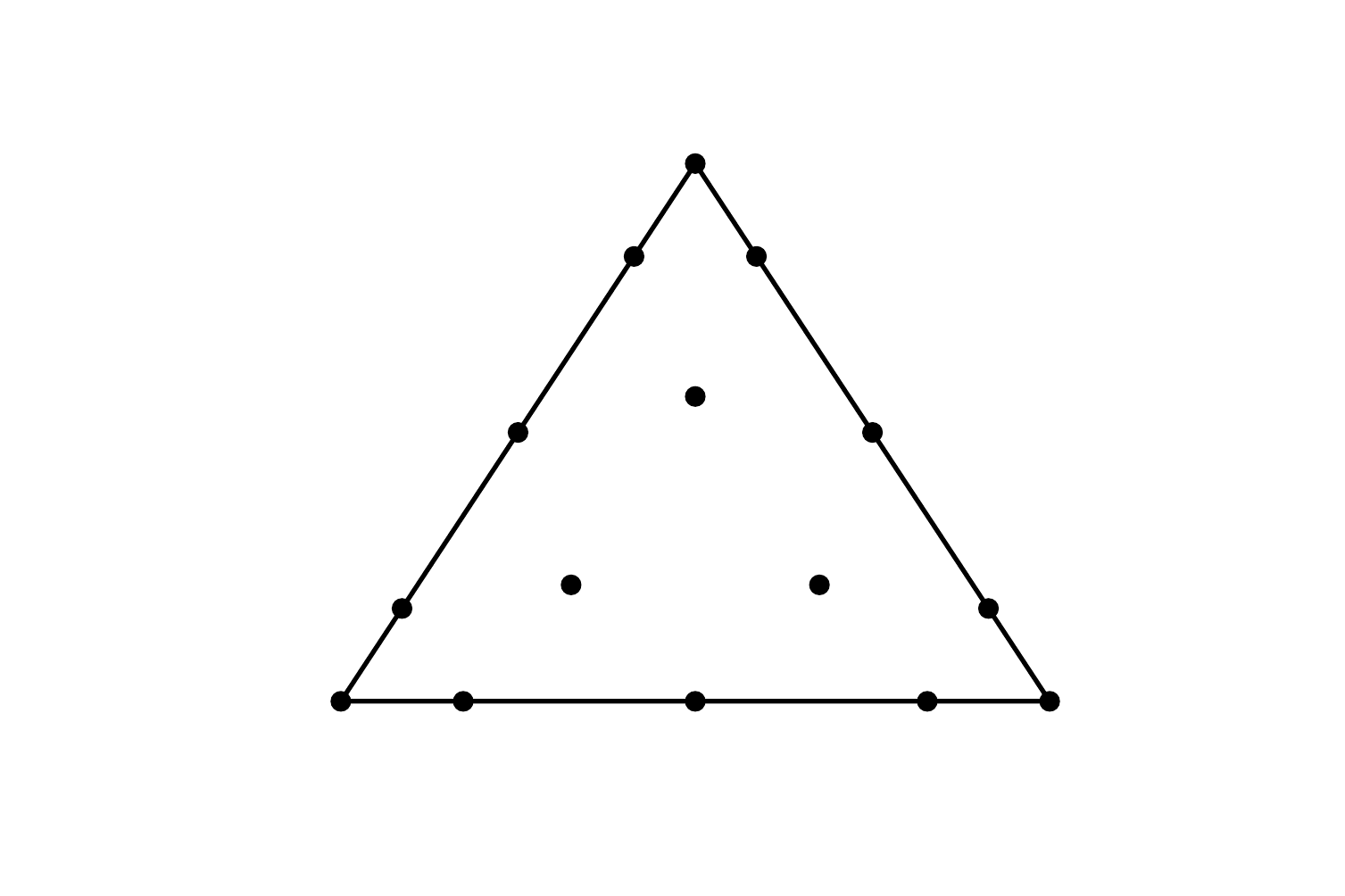}\hfill
 \caption{The 15 Points for Degree 4}
\end{figure}
\end{center}

It is shown in \cite{B83} that for $d=2,$ $F_4$ is a Fej\`er set.

\begin{prop} In dimension $d=2$ we have, for the set $F_4,$
\[\sum_{i=1}^N \ell_i^2({\bf x})\le 1,\,\,{\bf x}\in S_2.\]
\end{prop}

Extensive numerical calculations indicate that this is also the case in dimension $d=3.$

\begin{conj} In dimension $d=3$ we have,  for the set $F_4,$
\[\sum_{i=1}^N \ell_i^2({\bf x})\le 1,\,\,{\bf x}\in S_3.\]
\end{conj}

However,  it is {\it not} the case in dimension $d=4.$ Indeed,  by direct calculation we see that for ${\bf x}=(V_1+V_2+V_3+V_4+V_5)/5,$ the centroid of $S_4,$
\[\sum_{i=1}^N \ell_i^2({\bf x})=\frac{80243\sqrt{5}+9298842}{9453125}=1.00266...>1.\]
We nevertheless conjecture
\begin{conj} The sets $F_4$ are Fekete sets for all dimensions $d=1,2,\cdots.$
\end{conj}
As evidence for this conjecture we can prove that for dimensions $d=3,4,5,6,7$ the Lagrange polynomials for $F_4$ satisfy the necessary (but not sufficient) condition that they are all bounded by 1.  Indeed,  somewhat more is true.
\begin{prop} The Lagrange polynomials for the centroids of each 3-face,  
\[ \ell_{ijk\ell}=256\lambda_i\lambda_j\lambda_k\lambda_\ell
=256\lambda_i\lambda_j\lambda_k\lambda_\ell\]
are bounded by 1 in absolute value on $S_d$   for any dimension ($d\ge 3$).  Further,  the sum of all the {\bf other} Lagrange polynomials squared is bounded by 1 on $S_d$ for $d=3,4,5,6,7.$
\end{prop}
\noindent {\bf Proof}.  The first statement amounts to showing that
$256xyzw\le1$ for $x,y,z,w\ge0$ and $x+y+z+w\le1,$ but this is an elementary verification.

For the second statement,  we calculate,  just as in the proof of Prop. \ref{deg3case}, the sum of the stated Lagrange polynomials squared and again reduce to the sub-simplex with vertices $V_1,\cdots,V_d$ and the centroid
$\displaystyle{\frac{1}{d+1}\sum_{j=1}^{d+1}V_j.}$ Letting $\mu_j,$ $1\le j\le (d+1)$ denote the barycentric coordinates of this sub-simplex, we then have
\begin{align}\label{mus2}
\lambda_j&=\mu_j+\frac{1}{d+1}\mu_{d+1},\,\,1\le j\le d,\cr
\lambda_{d+1}&=\frac{1}{d+1}\mu_{d+1}.
\end{align}
Let $H({\bm \mu})$ be $1-K_4$ restricted to this sub-simplex,  i.e.,
\[H({\bm \mu})=\Bigl(\sum_{j=1}^{d+1}\lambda_j\Bigr)^8-K_4({\bm \lambda})\]
with ${\bm \lambda}$ given by (\ref{mus2}).  It is sufficient to prove that $H({\bm \mu})\ge 0$ for $\mu_j\ge0,$ $1\le j\le (d+1)$ for which it is, in turn, sufficient to show that
\[H({\bm \mu})-H(\widehat{{\bm \mu}},0)\ge 0,\quad \widehat{{\bm \mu}}:=(\mu_1,\mu_2,\cdots,\mu_d)\]
as the term $H(\widehat{{\bm \mu}},0)$ is the restriction of $H$ to the face $\mu_{d+1}=\lambda_{d+1}=0$ which we may assume to be positive by induction.  

For positivity it is again sufficient that, for some integer $r\ge0,$ all the coefficients of
\[ \Bigl(\sum_{j=1}^{d+1}\mu_j\Bigr)^r(H({\bm \mu})-H(\widehat{{\bm \mu}},0))\]
are non-negative.  This has been verified by means of computer algebra for degrees up to 7 with $r=11$ for dimension $d=3,$ $r=8$ for $d=4,$ $r=7$ for $d=5,6$ and $r=8$ for $d=7.$  $\square$

\begin{remark} One could of course continue this procedure dimension by dimension, but a general proof is much to be preferred.  $\square$
\end{remark}

Numerical computations indicate  that, at least for dimensions $d=3,4,5,$  the sum of the {\it fourth} powers  of the Lagrange polynomials is bounded by 1.

\begin{conj}  For all dimensions  $d\ge1,$ we have
\[\sum_{j=1}^N \ell_j^4({\bf x})\le 1,\,\,{\bf x}\in S_d\]
or, in other words,
\[ \|{\bm \ell}({\bf x})\|_4\le 1,\,\, {\bf x}\in S_d.\]
\end{conj}

\begin{remark} We use the fourth power in order to remain in the domain of polynomials.  However,  numerical computations indicate that there is a smaller exponent
\[r\approx 2.00217448\]
such that 
\[ \|{\bm \ell}({\bf x})\|_r\le 1,\,\, {\bf x}\in S_d,  \,\,d=4,5.\]
This leads us to formulate the following definition
\begin{definition} Suppose that $X\subset S_d$ is a subset of $N_n$ points for which $\|{\bm \ell}({\bf x})\|_\infty\le1,$ ${\bf x}\in S_d.$ We define the Fej\'er exponent of $X$ to be that $r\in [2,\infty]$ for which
\[\|{\bm \ell}({\bf x})\|_r\le 1,\,\, {\bf x}\in S_d\]
and for which for $r'<r,$ 
\[\max_{{\bf x}\in S_d} \|{\bm \ell}({\bf x})\|_{r'}>1.\]
\end{definition}
\noindent $\square$

We speculate that the Fekete points have minimal Fej\'er exponent.
\end{remark}

\section{Degree Five}

For degree 5 we are able to give explicit formulas for the points on  the  edges but can express the points in the interior of the 2-faces  only in terms of the roots of a certain polynomial.  Consequently this section will be largely numerical. 

We describe first the set of points $F_5$ in $\R^2.$ The dimension of the space of polynomials of degree at most five in two variables is $N_2={5+2\choose 2}=21.$ $F_5$ will consist of the 3 vertices together with 4 additional points on each of the 3 edges, leaving 21-15=6 points to be placed in the interior of the triangle.   The 4 interior edge points are just the interior univariate Fekete points of degree 5,  i.e.,  the zeros of $P_5'(x),$ the derivative of the Legendre polynomial of degree 5.  Now
\[P_5(x)= \frac{1}{8}\bigl( 63x^5 - 70x^3 + 15x\bigr)\]
and hence
\[P_5'(x)=\frac{15}{8}\bigl(21x^4-14x^2+1\bigr)\]
with zeros 
\[ z=\pm\sqrt{\frac{7\pm2\sqrt{7}}{21}}\]
with barycentric coordinates
\[\Bigr(\frac{1-z}{2},\frac{1+z}{2}\Bigr).\]

For the 6 interior points,  we assume that these are symmetric and have barycentric coordinates
\begin{align*}
&(u,u,1-2u),\,\,(u,1-2u,u),\,\,(1-2u,u,u)\cr
&(v,v,1-2v),\,\,(v,1-2v,v),\,\,(1-2v,v,v)
\end{align*}
for some parameters $u,v\in (0,1).$ Now,  the polynomials of degree 5 of the form $\lambda_1\lambda_2\lambda_3q({\bm \lambda}),$ for some quadrtic $q,$ are all zero at the boundary points.  Hence the Vandermonde determinant of all 21 points will decouple into a factor depending on the boundary points (given) and the Vandermonde determinant for the 6 interior points with basis $\lambda_1\lambda_2\lambda_3q({\bm \lambda})$  
for 
\[q\in \{ \lambda_1^2,\lambda_2^2,\lambda_3^2,\lambda_1\lambda_2,\lambda_1\lambda_3,\lambda_2\lambda_3\}.\]
(See e.g.  \cite{B91} for a discussion of this kind of factorization for Vandermonde determinants).  We obtain therefore that there is some constant $C$ such that
\begin{align*}
&{\rm vdm}=C(u^2(1-2u))^3(v^2(1-2v))^3\times \cr
&\left|\begin{array}{cccccc}
u^2&u^2&(1-2u)^2&u^2&u(1-2u)&u(1-2u)\cr
u^2&(1-2u)^2&u^2&u(1-2u)&u^2&u(1-2u)\cr
(1-2u)^2&u^2&u^2&u(1-2u)&u(1-2u)&u^2\cr
v^2&v^2&(1-2v)^2&v^2&v(1-2v)&v(1-2v)\cr
v^2&(1-2v)^2&v^2&v(1-2v)&v^2&v(1-2v)\cr
(1-2v)^2&v^2&v^2&v(1-2v)&v(1-2v)&v^2
\end{array}\right|\cr
&=C(u^2(1-2u))^3(v^2(1-2v))^3\times \cr
&(u-v)^3(3u-1)^2(3v-1)^2(3(u+v)-2)
\end{align*}
after a short calculation.

Hence we wish to maximize
\[p(u,v):=(u-v)^3u^6v^6(2u-1)^3(2v-1)^3(3u-1)^2(3v-1)^2(3(u+v)-2).\]
The partial derivatives are
\begin{align*}
\frac{\partial p}{\partial u}&=6(u-v)^2u^5v^6(2v-1)^3(3v-1)^2(2u-1)^2(3u-1)\times\cr
&\,\,\,\bigl(45u^4 + 6u^3v - 59u^3 - 33u^2v^2 + 17u^2v + 24u^2 + 21uv^2 - 13uv - 3u - 3v^2 + 2v\bigr)\cr
\frac{\partial p}{\partial v}&=6(u-v)^2u^6v^5(2u-1)^3(2v-1)^2(3u-1)^2(3v-1)\times \cr
&\,\,\,\bigl(33u^2v^2 - 21u^2v + 3u^2 - 6uv^3 - 17uv^2 + 13uv - 2u - 45v^4 + 59v^3 - 24v^2 + 3v\bigr).
\end{align*}
We seek to find the critical points given by the common zeros of the factors
\begin{align*}
q_1(u,v)&=45u^4 + 6u^3v - 59u^3 - 33u^2v^2 + 17u^2v + 24u^2 + 21uv^2 - 13uv - 3u - 3v^2 + 2v,\cr
q_2(u,v)&=33u^2v^2 - 21u^2v + 3u^2 - 6uv^3 - 17uv^2 + 13uv - 2u - 45v^4 + 59v^3 - 24v^2 + 3v.
\end{align*}
Upon calculating a Groebner basis of these polynomials one finds that these common zeros are determined by the univariate polynomial
\begin{align*}
q(v)&=2737800v^{10}- 8660340v^9 + 11981963v^8 - 9523289v^7\cr
&\,\,\, + 4800577v^6 - 1598001v^5 + 354286v^4 - 51415v^3 + 4649v^2 - 235v + 5
\end{align*}
whose roots may be calculated to any precision desired. By trial and error we find that the roots near
\begin{align*}
u&=0.148019471315134,\cr
v&=0.420825539292557
\end{align*}
give the largest determinant, and we use these.

\begin{center}\begin{figure}[t!]\centering
 \includegraphics[width=.6\textwidth,height=.3\textheight]{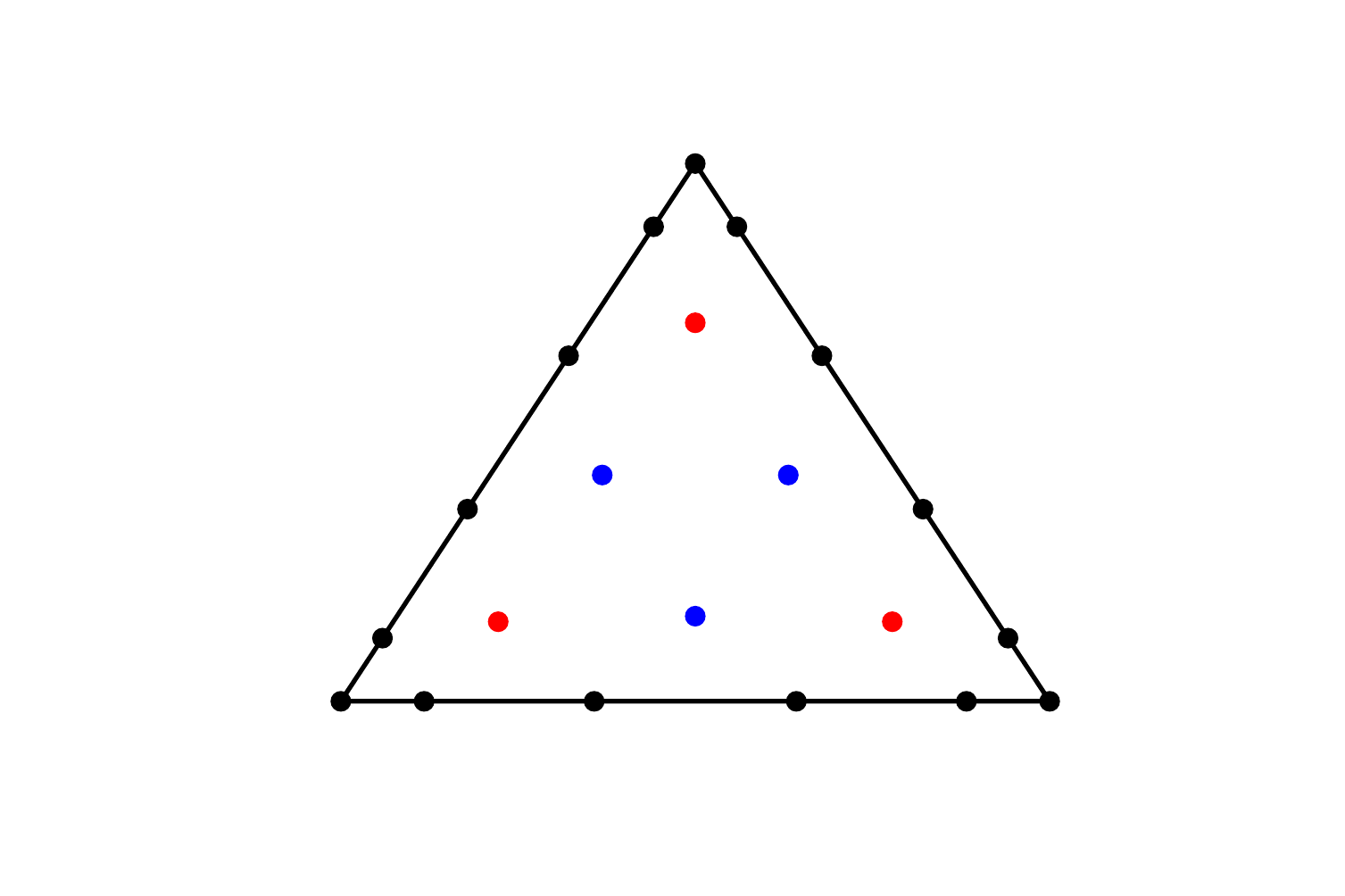}\hfill
 \caption{The 21 Points for Degree 5}
\end{figure}
\end{center}

\begin{remark}
Although the 6 interior points appear visually to all be on an interior {\it triangle}, close inspection of the coordinates reveals that this is not the case. $\square$
\end{remark}

\begin{example} For dimension $d=2,$ the set $F_5$ of the 21 points defined above do {\bf not} form a Fej\'er set. Indeed,
\[\max_{{\bf x}\in S_2}\sum_{i=1}^{21} \ell_i^2({\bf x})\approx 1.2246\]
which is attained at the centroid $(1/3,1/3,1/3).$

However,  we do believe that they are Fekete points.  Indeed it's Fej\'er exponent is
\[r\approx 2.2513.\]
In particular,  all the Lagrange polynomials are bounded by 1 in absolute value on $S_2.$ $\square$
\end{example}

We now proceed to dimension $d=3.$ Here the dimension of the space of polynomials of degree at most 5 is
\[N_5={5+3\choose 3}=56.\]
The set $F_5$ is such that restricted to any 2-face we obtain its two-dimensional version,  i.e.,  we place 1 point at each of the 4 vertices,  4 interior points on each of the 6 edges and 6 interior points on each of the 4 2-faces,  for a total of $4+4\times 6+6\times 4=52$ boundary points.  This leaves $56-52=4$ points to be placed in the interior of the simplex.  By symmetry we assume that these are 
\begin{align}\label{intpts5}
&(1-3w,w,w,w)\cr
&(w,1-3w,w,w),\cr
&(w,w,1-3w,w),\cr
&(w,w,w,1-3w)
\end{align}
for some parameter $w\in (0,1/3).$

The polynomials of the form 
\[\lambda_1\lambda_2\lambda_3\lambda_4 q({\bm \lambda}),\,\, {\rm deg}(q)\le 1\]
are all zero on the boundary.  Hence (cf.  \cite{B91}) the Vandermonde determinant for the 56 points will factor into a constant (depending on the boundary points) times the Vandermonde determinant for the 4 points (\ref{intpts5}) with basis
\[\lambda_1\lambda_2\lambda_3\lambda_4 \{\lambda_1,\lambda_2,\lambda_3,\lambda_4\}\]
i.e.,
\begin{align*}
&(w^3(1-3w))^4\,\,\left|\begin{array}{cccc}
1-3w&w&w&w\cr
w&1-3w&w&w\cr
w&w&1-3w&w\cr
w&w&w&1-3w \end{array}\right|\cr
&=(w^3(1-3w))^4(1-4w)^3.
\end{align*}
The derivative of this expression with respect to $w$
is
\[ 12w^{11}(1-3w)^3(1-4w)^2(19w^2 - 9w + 1)\]
with non-extraneous crtical points
\[ w=\frac{9\pm\sqrt{5}}{38}.\]
Now, it turns out that the critical point $(9+\sqrt{5})/38$ is a {\it local} maximum for the absolute value of the determinant and $(9-\sqrt{5})/38$ the {\it global} maximum.  Hence we take
\begin{equation}\label{w}
w=\frac{9-\sqrt{5}}{38}.
\end{equation}

\begin{example} There is a curious behaviour here.  For {\it both} choises of $w=(9\pm\sqrt{5})/38,$ the Lagrange polynomials are all bounded by one in absoliute value,  i.e.,
\[\|{\bm \ell}({\bf x})\|_\infty \le 1,\,\,{\bf x}\in S_3.\]
However,  the choice of $w=(9-\sqrt{5})/38$ produces a strictly larger Vandermonde determinant and hence only with this choice of $w$ can 
$F_5$ be a Fekete set,  which we conjecture to be the case.  

It is also interesting to note that with $w=(9+\sqrt{5})/38, $ the Fej\'er exponent is $>4,$ wheras for  $w=(9-\sqrt{5})/38, $ the Fej\'er exponent remains, as in dimension $d=2,$
\[r\approx 2.2513.\]
$\square$
\end{example}

Moving to dimension $d=4,$ 
\[N_5={5+4\choose 4}=126.\]
We again place points on the boundary so that restricted to any 3-face we have $F_5$ for dimension 3.  In this way we have 5 vertex points,  4 interior points on each of the ${5\choose 2}=10$ edges,  6 interior points on each of the ${5 \choose 3}=10$ 2-faces and 4 interior points on each of the ${5 \choose 4}=5$ 3-faces,  for  a total of $5+4\times 10+6\times 10+4\times 5=125$ boundary points so that there is but $1=126-125$ to be placed in the interior of the simplex.  We put this point at the centroid $(1/5,1/5,1/5,1/5,1/5).$

\begin{example} We conjecture that the set $F_5$ of 126 points so constructed is a Fekete set.  Indeed it remains the case that the Fej\'er exponent is
\[r\approx 2.2513\]
and so, in particular,  all the Lagrange polynomials are bounded by 1 in absolute value.  $\square$
\end{example}

In general we let $F_5$ be the set of points  consisting of
\begin{itemize}
\item the ${d+1 \choose 5}$ centres of each  four dimesnional face: \[\frac{1}{5}(V_i+V_j+V_k+V_\ell+V_m),\,\, i<j<k<\ell<m\]
\item the $4{d+1 \choose  4}$ vertices of the simplex in the interior of each three dimensional face:
\begin{align*}
&(1-3w)V_i+wV_j+wV_k+wV_\ell,\cr
&wV_i+(1-3w)V_j+wV_k+wV_\ell,\cr
&wV_i+wV_j+(1-3w)V_k+wV_\ell,\cr
&wV_i+wV_j+wV_k+(1-3w)V_\ell
\end{align*}
for $i<j<k<\ell$ and $w=(9-\sqrt{5})/38.$
\item the $6{d+1 \choose 3}$ points in each 2-face:
\begin{align*}
&(1-2u)V_i+uV_j+uV_k,\cr
&uV_i+(1-2u)V_j+uV_k,\cr
&uV_i+uV_j+(1-2u)V_k,\cr
&(1-2v)V_i+vV_j+vV_k,\cr
&vV_i+(1-2v)V_j+vV_k,\cr
&vV_i+vV_j+(1-2v)V_k
\end{align*}
for $i<j<k$ and $u,v$ two of the roots of
\begin{align*}
q(v)&=2737800v^{10}- 8660340v^9 + 11981963v^8 - 9523289v^7\cr
&\,\,\, + 4800577v^6 - 1598001v^5 + 354286v^4 - 51415v^3 + 4649v^2 - 235v + 5
\end{align*}
near $u=0.148019471315134,$ 
$v=0.420825539292557.$
\item the $4{d+1 \choose 2}$ interior edge points:
\[\lambda V_i+(1-\lambda)V_j\]
for $i<j$ and $\lambda = (1-z)/2,$ $\displaystyle{z=\pm\sqrt{\frac{7\pm2\sqrt{7}}{21}}}$
\item the $(d+1)$ vertices $V_i.$
\end{itemize}
We note that the cardinality of $F_5$ is
\begin{align*}
\#(F_5)&={d+1 \choose 5}+4{d+1 \choose 4}+6{d+1 \choose 3}+4{d+1 \choose 2}+{d+1 \choose 1}\cr
&={d+1\choose 5}\cr
&=N_5.
\end{align*}

\section*{Appendix A}

\begin{lstlisting}
%
% test the positivity of 1-K_3 on the d-dimensional sub-simplex 
% with vertices V_1,...V_d and the barycentre (V_1+...V_{d+1})/(d+1)
% of the simplex with vertices V_1,...V_{d+1} in R^d
%
clear all
syms K t s1 s2
d=input('Desired dimension: ')
x=sym('x',[d+1 1]);   % barycentric coordinates for the whole simplex
y=sym('y',[d+1 1]);   % barycentric coordinates for the sub-simplex
%
% find K, the sum of the squares of the Lagrange polynomials on the simplex
%
t=(1+1/sym(sqrt(5)))/2;
s1=sum(x);
s2=sum(x.^2);
K=0;
for i=1:(d+1)
    K=K+(x(i)/2*(12*x(i)^2-12*x(i)*s1+3*s1^2-s2))^2;
end
for i=1:(d+1)
    for j=1:(i-1)
        K=K+(5*x(i)*x(j)*((1+sym(sqrt(5))*t)*x(i)+(2-sym(sqrt(5))*t) ...
            *x(j)-s1))^2;
        K=K+(5*x(i)*x(j)*((1+sym(sqrt(5))*t)*x(j)+(2-sym(sqrt(5))*t) ...
            *x(i)-s1))^2;
    end
end
for i=1:(d+1)
    for j=1:(i-1)
        for k=1:(j-1)
            K=K+(27*x(i)*x(j)*x(k))^2;
        end
    end
end
K=expand(K);
f=expand(s1^6-K); % want this to be positive
%
% change to barycentric coordinates on the sub-simplex
% x_j=y_j+y_{d+1}/(d+1), j=1...d, x_{d+1}=y_{d+1}/(d+1)
%
g=f;
for j=1:d, g=subs(g,x(j),y(j)+y(d+1)/(d+1)); end
g=expand(subs(g,x(d+1),y(d+1)/(d+1)));
g1=g;
% 
% now subtract off the part corresponding to the face y_{d+1}=x_{d+1}=0
% which can be assumed to be positive by induction
%
g=expand(g-subs(g,y(d+1),0));
%
% we now test that this remainder is positive on the sub-simplex
% by raising its degree (multiplying by the sum of the coordinates to a
% power, and verifying that all the coefficents are non-negative
%
r=input('Raise the degree by: ')
h=expand((sum(y))^r*g);
mm=min(coeffs(h));    % if this is positive all coefficients are non-negative
disp('Minimum (non-zero) coefficient')
mm
\end{lstlisting}

\bigskip

\section*{Acknowledgements}
 RITA 
``Research ITalian network on Approximation''.

\end{document}